% plaintex
% public copy

\magnification=\magstep1

\overfullrule=0pt
\parskip=3pt
\baselineskip=14.8pt
\abovedisplayskip=5pt
\belowdisplayskip=5pt

\font\bbigfont=cmbx12

\def\today{\ifcase\month\or January\or February\or March\or
	April\or May\or June\or July\or August\or September\or
	October\or November\or December\fi
	\space\number\day, \number\year}

\def \reg{\hbox{reg$\,$}}
\def \qed{\qquad \vrule width 1.6mm height 1.6mm depth 0mm}

\def \sumi{{\sum_{j=1}^m I_j}}
\def \sumin{{\sum_{j=1}^m I_j^n}}
\def \sumip{{\bigl(\sum_{j=1}^m I_j\bigr)}}
\def \sumipn{{\bigl(\sum_{j=1}^m I_j^n\bigr)}}

\newcount\sectno \sectno=0
\newcount\thmno \thmno=0
\def \section#1{\bigskip\bigskip
	\global\advance\sectno by 1 \global\thmno=0
	\noindent{\bf \the\sectno. #1}}
\def \thmline#1{\vskip 5pt
	\global\advance\thmno by 1
	\noindent{\bf (\the\sectno.\the\thmno)\ #1:}\ \ %
	\bgroup \it}
\def \endthm{\egroup \vskip 5pt}

\def \thm{\thmline{Theorem}}
\def \cor{\thmline{Corollary}}
\def \lemma{\thmline{Lemma}}

\def \proof{\vskip 0mm\noindent {\sl Proof:\ \ }}

\def\label#1{\global\edef#1{\the\sectno.\the\thmno}}
\def\eqlabel#1{\global\advance\eqnum by1\global\edef#1{\the\sectno.\the\eqnum}}

\smallskip
\centerline{\bbigfont
Linear Bounds on Growth of Associated Primes}
%for monomial Ideals }
\bigskip
\centerline{Karen E.\ Smith}
\centerline{Department of Mathematics}
\centerline{Massachusetts Institute of Technology 2-167}
\centerline{Cambridge, MA 02139}
\centerline{{\tt kesmith@math.mit.edu}}
\medskip
\centerline{and}
\medskip
\centerline{Irena Swanson}
\centerline{Department of Mathematical Sciences}
\centerline{New Mexico State University}
\centerline{Las Cruces, NM 88003-8001}
\centerline{{\tt iswanson@nmsu.edu}}
\bigskip
%\centerline{\today}
%\bigskip

Let $I$ be an ideal in a commutative Noetherian ring $R$.
We study the growth of associated primary components of powers of $I$.
Swanson [S] and Heinzer and Swanson [HS]
have shown that there exists an integer $l$ such that
each power of $I$ has a primary decomposition
$$
I^n = q_1 \cap \cdots \cap q_s
$$
where $(\sqrt q_i)^{nl} \subset q_i$ for each of the primary components
$q_i$. In other words, the primary components of $I^n$
``grow linearly'' in $n$. An interesting question of practical
concern is {\it what is \/} $l$?

Swanson was partially motivated by the analogous question for
``Frobenius powers'',
which has important applications to the theory of tight closure
(the so-called ``localization problem;'' see [HH, Section 4]).
Suppose now that $R$ has prime characteristic $p>0$
and $q$ is a power of $p$.
The Frobenius power $I^{[q]}$ of an ideal $I$
is the ideal generated by all the $q$th
powers of the elements (equivalently, the generators) of $I$.

We ask the same question: can we find an integer $l$ such that for each $q$
there is a primary decomposition of $I^{[q]}$
$$
I^{[q]} = q_1 \cap \cdots \cap q_s
$$
with each $(\sqrt q_i)^{[q]l} \subset q_i$?
This Frobenius analog of [S] may be quite difficult.
One difficulty is that, unlike with the case of ordinary powers of $I$
(proved in [R]),
the Frobenius powers of $I$ can have infinitely many distinct associated primes
as $q$ varies over all powers of $p$.
An ideal with this property was first found by Katzman [K].
We analyze the primary decompositions of Frobenius powers
of Katzman's ideal in Section~2
and prove that the question above has an affirmative answer for this ideal.

In Section~1 we solve the problem of linear growth
for primary components of both ordinary and Frobenius powers
of monomial ideals in polynomial rings modulo a monomial ideal
%a particular class of ideals in a particular class of rings by
by giving an explicit value $l$ that works.
More precisely,
we analyze the ideals generated by monomials in a ring of the form $R = S/J$,
where $S$ is a polynomial ring over a field
and $J$ is an ideal generated by monomials in the given variables.
Because primary decompositions of monomial ideals are easy to understand,
we are able to give a constructive bound that works in both
the ordinary and the Frobenius power case.
We point out how this can be used to verify that tight closure 
commutes with localization in this special case.

In Section 3 we discuss the Castelnuovo-Mumford regularity of
powers of monomial ideals.
For a given monomial ideal $I$
we explicitly find an integer $B$
such that the regularity of the $n$th power of $I$
is bounded above by $Bn$.
We believe that in general $B$ can be improved.
However, in many cases the given $B$ is a sharp upper bound.

%In Section 2, we examine Katzman's example of an ideal whose
%Frobenius powers have infinitely many associated primes.
%Even though his
%example shows that there can be new associated primes for $I^{[q]}$
%as $q$ grows,
%we show that nevertheless,
%the linear growth property
%discussed at the beginning of this paper
%is still valid in his example.

\section{\bbigfont Primary decompositions and localization of tight closure
on monomial ideals}

Throughout in this section
$S$ is a polynomial ring $k[x_1, \ldots, x_d]$,
where $k$ is a field and $x_1, \ldots, x_d$ are variables over $k$.
Let $I$ and $J$ be ideals generated by monomials in these variables.

First some notation.
An ideal is said to be {\it irreducible} if it
cannot be written as an intersection of two strictly
larger ideals.
It is well-known and easy to show that every irreducible monomial ideal
is of the form
$(x_{i_1}^{a_1}, x_{i_2}^{a_2}, \dots, x_{i_t}^{a_t})$
for some integers $t \le d$,
$i_1 < i_2 < \cdots < i_t$,
and $a_i \ge 1$.
The irreducible monomial ideals can be indexed by monomials:
$(x_{i_1}^{a_1}, x_{i_2}^{a_2}, \dots, x_{i_t}^{a_t})$
corresponds uniquely to the monomial
$x_{i_1}^{a_1}x_{i_2}^{a_2} \cdots x_{i_t}^{a_t}$,
where each $a_i$ is a positive integer.
This is a one-to-one correspondence,
so we will write an irreducible monomial ideal as $J_m$ for a monomial $m$.

Consider an arbitrary monomial ideal $I$
in the polynomial ring $k[x_1, x_2, \dots, x_d]$.
The ideal $I$ can be decomposed into irreducible ideals
$$
I = J_{m_1} \cap J_{m_2} \cap \cdots \cap J_{m_r}
$$
for some monomials $m_1, \dots, m_r$
(see for example Eagon-Hochster [EH], Sturmfels-Trung-Vogel [STV],
or Heinzer-Ratliff-Shah [HRS]).
Assume this decomposition is {\it minimal,}
so that no $J_{m_i}$ can be omitted,
and no $m_i$ can be replaced with a proper monomial factor.

\lemma
With irreducible decomposition of $I$ as above,
let $l$ be the largest of the powers of a variable appearing
as minimal generators of the $J_{m_i}$ in this minimal decomposition.
(Equivalently,
$l$ is the largest exponent of any variable appearing in the indexing monomials
$m_1, m_2, \dots, m_r$.)
Let $k$ be the largest exponent of a variable in a minimal
generating set of $I$.
Then $l$ equals $k$.
\endthm

\proof
Note that
$I$ is generated by least common multiples
of elements $\alpha_1, \ldots, \alpha_r$
as the $\alpha_i$ run through the monomial generators of $J_{m_i}$.
Thus in the monomial generating set of $I$ obtained in this way,
all the variables have exponent at most $l$,
so each minimal generator of $I$ will have all exponents bounded above by $l$.
This proves that $k \le l$.

Suppose that $k < l$.
By reindexing, if necessary,
$l$ may be assumed to be the exponent of ${x_1}$ in $m_1$.
Let $m_1' = m_1/x_1^l$.
Then certainly
$$
I = J_{m_1} \cap J_{m_2} \cap \cdots \cap J_{m_r} \ \ \supseteq \ \ I' =
J_{m_1'} \cap J_{m_2} \cap \cdots \cap J_{m_r}.
$$
%Note that
%the associated prime of $J_{m_1'}$ is strictly contained in
%the associated prime of $J_{m_1}$,
%so that the inclusion $I' \subseteq I$ is proper.
%However, we now show that $I' = I$:
By assumption,
for each minimal generator $m = x_1^{a_1}\cdots x_d^{a_d}$ of $I$,
we have $a_1 < l$.
Thus $m$ is a multiple of the monomials in $J_{m_1}$ other than $x_1^l$
which means that
$m$ lies in $J_{m_1'}$.
This says that each of the minimal generators for $I$ is actually in $I'$,
whence $I = I'$,
contradicting the minimality of our irreducible decomposition.
\qed

This lemma makes the proof of the linear growth of primary decompositions
of both ordinary and Frobenius powers very easy:
for let $I$ and $J$ be monomial ideals in $S$.
Let $l$ be the largest exponent appearing
in a set of minimal monomial generators for $I+J$.
Then the corresponding largest exponent for
$I^q + J$ and $I^{[q]} + J$ is no more than $ql$.
% in the latter case $q$ is interpreted to be a power of the characteristic $p$
%base field).
If we decompose $I^q + J$ (or $ I^{[q]} + J$)
efficiently into irreducibles,
then each irreducible component $J_m$ involves powers of variables of degree
at most $ql$.
Clearly $(\sqrt{J_m})^{qld} \subseteq J_m$.
Thus:

\thm
\label{\thmFrob}
Let $I$ and $J$ be monomial ideals in $S = k[x_1, \ldots, x_d]$.
Let $l$ be the largest power of any variable occurring
in a generating set of $I$ or of $J$.
Then for each $n$ there exists a primary decomposition
$I^n + J = q_1 \cap \cdots \cap q_s$
such that
$\sqrt{q_i}^{nld} \subseteq q_i$ for all $i$.
If the characteristic $p$ of $k$ is positive,
then for any power $q$ of $p$
there exists a primary decomposition
$I^{[q]} + J = q_1 \cap \cdots \cap q_s$
such that
$\sqrt{q_i}^{nld} \subseteq q_i$ for all $i$.
Moreover,
for $n$ and $q$ sufficiently large
we may take $l$ to be 
the largest power of any variable occurring in the minimal generators for $I$.
\endthm

\proof
Note that a primary decomposition is obtained from an irreducible one
by simply intersecting the irreducible components with the same radical.
The result is now immediate from the lemma and the discussion above.
\qed

This proves the linear growth of primary decompositions
of ordinary and Frobenius powers of a monomial ideal $I$
in $S/J$,
where $J$ is a monomial ideal in the polynomial ring
$S = k[x_1, \ldots, x_d]$ over a field $k$.

This theorem also has an immediate application to the localization
problem in tight closure.
Tight closure is a closure operation performed
on ideals in a commutative Noetherian ring of prime characteristic.
For an ideal $I \in R$, an element $z \in R$ is in the tight closure
$I^*$ of $I$ if there exists an element $c$ not in any minimal prime of $R$
such that $cz^q \in I^{[q]}$ for all $q = p^n >> 0$. See [HH].

One of the most persistent open problems in the theory of tight closure
has been the ``localization problem:'' given an ideal $I\subset R$ and
a multiplicative system $U \subset R$, is $I^*U^{-1}R = (IU^{-1}R)^*$?
The difficulty in understanding the growth of associated primes
of monomial ideals is one of the key obstructions to 
settling this problem.In the next corollary, we indicate how our 
result on linear growth of associated primes can be applied to 
see that tight closure commutes with localization for monomial rings.

\cor
\label{\corlocalization}
Let $I$ be an ideal generated by monomials in $R = S/J$ where $J$ is a monomial
ideal in $S = k[x_1, \dots, x_d]$. Then for any $u \in S$,
$$I^*R_u = (IR_u)^*;$$
that is, tight closure commutes with localization
at multiplicatively closed sets of the form $\{1, u, u^2, \ldots \}$
for monomial ideals
in rings of the type $S/J$.
\endthm

\proof
The inclusion $I^*(R_u) \subset (IR_u)^*$ is obvious.
Suppose that ${z \over 1} \in (IR_u)^*$.
Then for each $q = p^n$, there exists an integer $N(n)$
such that $cz^qu^{N(n)} \in I^{[q]}$.
{\it A priori,\/} we know
only that $c$ is not in any minimal prime of $R$ not containing $u$,
but by replacing $c$ by $c + \delta$ where $\delta$ belongs precisely
to those minimal primes not containing $c$ and no others, we may assume
that $c$ is in no minimal prime of $R$ (see the proof of Proposition 4.14
in [HH]).
If $cz^q$ is in each primary component of $I^{[q]}$,
then $cz^q \in I^{[q]}$, and we are done.
Otherwise,
the element $u$ is
in every associated prime of $I^{[q]}$
such that $cz^q$ is not in the corresponding primary component.
By Theorem \thmFrob,
$u^{qN}$ is in each of these primary components,
where $N$ is a fixed integer independent of $q$.
Therefore, $c(zu^N)^q \in I^{[q]}$,
whence $zu^N \in I^*$ and $z \in I^*R_u$.
\qed

Despite the interesting method of proof, which we hope may eventually be 
more broadly applicable, the  corollary
does not really prove anything new about localization.
It is easy to see that tight closure  commutes with localization
at an arbitrary multiplicative system  for any ideal $I$ 
in a monomial ring $R$: in
these rings, tight closure has a simple description as 
 $I^* = \cap_{P \in minspec R} (I + P)$, and localization is immediate.
This description of tight closure 
 follows from  the fact that tight closure can be computed 
modulo minimal primes and  the fact 
that all ideals are tightly closed in a regular ring.
More explicitly, 
 $z \in I^* $ if and only if the image of $z$ is in $(IR/P)^*$ for
 each minimal prime of $R$, and  when $R$ is a monomial ring, each $R/P$ is
a polynomial ring, so $(IR/P)^* = IR/P$.
 See also [T], [K2].

%\vfill\eject
\section{\bbigfont Primary decompositions of Katzman's ideal}

Let $k$ be a field of characteristic $p > 0$,
$t$, $x$ and $y$ indeterminates over $k$,
$A = k[t]$ and $R = A[x,y]$.
Katzman's example is as follows.
Set $I_q = (x^q, y^q, xy(x-y)(x-ty))R$.
As $q$ varies over powers of $p$ (or all integers),
the set of associated primes of the $I_q$ is infinite.
In particular,
this means that the set of associated primes
of all the Frobenius powers of $I = (x, y)$
in the ring $k[t, x, y]/(xy(x-y)(x-ty)$
is infinite.
%as $q$ ranges through all powers of $p$.

We now show that although the set of associated primes is infinite,
there nonetheless exists an integer $l$ such that for each $q$,
there is a primary decomposition of $I^{[q]}$:
$$
I^{[q]} = q_1 \cap \cdots \cap q_s
$$
with $(\sqrt q_i)^{lq} \subset q_i$ for all $q_i$.
In fact,
we show that $l = 2$.

\medskip
\def \gq{x^2y^{q-1}}
Consider the elements
 $\tau_q = 1 + t + t^2 + \cdots + t^{q-2}$ and
$G_q = \gq$ of $R$.
Katzman showed that
$\tau_q G_q \in I_q$.
Define $J_q = I_q + \gq R$.
We show below that $J_q$ 
is a primary component of $I_q$.

\lemma
\label{\lemmaKatzman}
Let $f$ be any nonzero element in $A$.
Then $J_q : f = J_q$ for all $q$.
\endthm

\proof
We have to prove that
$J_q : f \subseteq J_q$.
Let
$\alpha \in J_q : f$.
Without loss of generality we may assume that $\alpha$ is homogeneous of
degree $d$ under the $(x,y)$-grading.
Note that $(x,y)^{q+1} \subseteq J_q$.
%(x^q, y^q) + (x,y)^{q+1} + (xy(x-y)(x-ty))$,
Thus without loss of generality we may assume that
$d \le q$ and $f \alpha$ is a multiple of $(xy(x-y)(x-ty))$.
But $f, xy(x-y)(x-ty)$ is a regular sequence,
so $\alpha \in (xy(x-y)(x-ty))$.
\qed

In particular, $J_q : \tau_q = J_q$,
hence $I_q : \tau_q \subseteq J_q : \tau_q = J_q \subseteq I_q : \tau_q$.
Moreover,
$I_q : \tau_q^2 = (I_q : \tau_q) : \tau_q = J_q : \tau_q = J_q = I_q : \tau_q$,
which implies that
$I_q = (I_q : \tau_q) \cap (I_q + (\tau_q))$.

So in order to analyze the primary components of $I_q$
it suffices to analyze the primary components of
$I_q : \tau_q$ and $I_q + (\tau_q)$ separately.

First we %HERE I added we
analyze the primary components of $J_q = I_q : \tau_q$.
As $J_q$ is homogeneous under the $(x,y)$-grading,
all the prime ideals associated to it are also homogeneous.
We claim that $J_q$ is $(x,y)$-primary.
If it is not,
then there exists an element in $A$ which is a zero-divisor on $R/J_q$.
But this is impossible by Lemma~\lemmaKatzman.
Now observe that $(x,y)^{q+1}$ and hence $(x,y)^{2q}$
is contained in $J_q$.
This means that the
linear growth property holds for the primary component $J_q$ of $I_q$;
namely $(\sqrt {J_q})^{2q} \subset J_q$.

Now consider the primary components of $I_q + (\tau_q)$.
Because the radical of this ideal includes $x, y$ and $\tau_q$,
it is clear that this ideal is height three in $k[t, x, y]$, and thus cannot
have any embedded primary components.
Let $\tau_q = \prod_{i=1}^{r} \sigma_i$
be a factorization of $\tau_q \in k[t]$ into distinct irreducible polynomials.
%HERE I added distinct above and changed below a little
Each $\sqrt{I_q + (\sigma_i)} = (x, y, \sigma_i)$
is a maximal ideal containing $I_q + (\tau_q)$
and these are the only maximal ideals containing 
$I_q + (\tau_q)$.
An $(x,y,\sigma_i)$-primary component of $I_q + (\tau_q)$
has to contain $\sigma_i$,
so we get a primary decomposition
$$
I_q + (\tau_q) = (I_q + (\sigma_1)) \cap (I_q + (\sigma_2)) \cap \cdots \cap
(I_q + (\sigma_r)).
$$
Now it is clear that $(x, y, \sigma_i)^{2q} \subset (x^q, y^q, \sigma_i)
\subset I_q + (\sigma_i)$.
Therefore the same linear bound $2q$
that worked for the primary component $J_q$
also works for $I_q + (\tau_q)$.

In summary, Katzman's examples
$I_q$ decompose as
$$
I_q = J_q \cap q_1 \cap \cdots \cap q_s
$$
where $J_q = I_q: \tau_q = I_q + (x^2y^{q-2})$
is primary with $(\sqrt {J_q})^{2q} \subset J_q$, and each
$q_i = I_q + (\sigma_i)$ is primary with $(\sqrt{q_i})^{2q} \subset q_i$.

%\vfill\eject
\section{\bbigfont Castelnuovo-Mumford regularity of monomial ideals}

This section is a preliminary attempt at understanding how
for a given monomial ideal $I$,
the Castelnuovo-Mumford regularity varies with powers of $I$.
If $I$ is Borel-fixed,
then the Eliahou-Kervaire resolution gives that
$\reg(I^n) \le n \reg(I)$ (see [EK]).
We do not know whether this is true for arbitrary monomial ideals.

Here we do the following:
given a monomial ideal $I$ in $k[x_1, \ldots, x_d]$,
we want to find the integer $B$
such that for all $n$,
the Castelnuovo-Mumford regularity of $I^n$ is bounded above by $Bn$.
Such an integer $B$ exists by Theorem~3.6 in [S]
but the arguments in [S] do not show how to calculate $B$.
We show in this section that for monomial ideals 
one can calculate such a $B$.

In the course of the proof
we found it necessary to determine an upper bound on the
Castelnuovo-Mumford regularity for a more general class of ideals.
Namely we prove:

\thm
\label{\thmreg}
Let $R$ be $k[x_1, \ldots, x_d]$,
a polynomial ring in $d$ variables over a field $k$.
Let $I_1, \ldots, I_m$ be monomial ideals in $R$.
Let $l$ be the largest exponent of a variable occurring in
any of the generating sets for the $I_j$.
Also assume $x_1, \ldots, x_r$ all lie in the radical of $\sum_j I_j$.
For a subset $S$ of the variables and $x_q \in S$,
define
$I_{S,x_q} = \sum _{x_i \in S\setminus \{x_q\}} x_i^lR$.
Let
$$
L = \max_{S \subseteq \{x_{r+1}, \ldots, x_d\}}
\max_{x_q \in S}
\{
(d-|S|-r)l + (m+|S|-2+d) \reg \Bigl(
(\sumip :x_q^l) + I_{S,x_q} \Bigr)
\}.
$$
Then $\reg \Bigl( \sum_{j=1}^m I_j^n \Bigr) \le n \max\{dl, L\}$.
\endthm

Note that $dl$ and $L$ are both computable
and that
the ideals $(\sumip :x_q^l) + I_{S,x_q}$
involve at most $d - 1$ variables.

Thus for a single monomial ideal $I$
such that $l$ is the largest exponent of a variable occurring
in a monomial generating set
and such that $x_1, \ldots, x_r$ all lie in the radical of $I$
we get that
$$
reg(I^n) \le
n
\max \{dl,
\max_{S \subseteq \{x_{r+1}, \ldots, x_d\}}
\max_{x_q \in S}
\{
(d-|S|-r)l + (|S|-1+d) \reg \Bigl(
(I :x_q^l) + I_{S,x_q} \Bigr)
\}.
$$

Before we prove Theorem~\thmreg\ we need a few lemmas:

\lemma
\label{\lmeisen}
Let $0 \longrightarrow A
\longrightarrow B
\longrightarrow C
\longrightarrow 0$
be a short exact sequence of graded finitely generated $R$-modules
such that all the maps are homogeneous of degree 0.
Then
\itemitem{(i)}
$\reg A \le \max\{ \reg B, \reg C + 1\}$,
\itemitem{(ii)}
$\reg B \le \max\{ \reg A, \reg C \}$.
\endthm

\proof
See for example Corollary 20.19 in Eisenbud [E].
\qed

\lemma
\label{\lmlm}
Let $I$ be a homogeneous ideal in $R$,
$x_q$ a variable and $l$ a positive integer.
Then
\itemitem{(i)}
$\reg(I \cap x_q^lR) = l + \reg (I : x_q^l)$.
\itemitem{(ii)}
$\reg(I : x_q^l) \le \max \{ 0, \reg(I) - l, \reg(I + x_q^lR) + 1 - l\}$.
\itemitem{(iii)}
$\reg(I) \le \max \{ \reg(I : x_q^l) + l, \reg(I + x_q^lR)\}$.
\endthm

\proof
(i) follows from the elementary facts
$I \cap x_q^lR = x_q^l (I : x_q^l)$
and
$\reg (x_q^l J) = l + \reg J$
for any ideal $J$.

For (ii) and (iii),
we use
Lemma~\lmeisen, (i) and the short exact sequence
$$
0 \longrightarrow I \cap x_q^l R
\longrightarrow I \oplus x_q^l R
\longrightarrow I + x_q^l R
\longrightarrow 0,
$$
to get that
$$
\eqalign{
\reg(I:x_q^l)
&=
\reg(I \cap x_q^lR) - l \cr
&\le
\max \{ \reg (I \oplus x_q^l R) - l, \reg(I + x_q^lR) + 1 - l\} \cr
&=
\max \{ \reg (I) - l, \reg(x_q^l) - l, 
\reg(I + x_q^lR) + 1 - l\},
\cr
}
$$
which proves (ii).
Similarly, (iii) follows.  \qed

Now we prove the main theorem of this section,
namely Theorem~\thmreg:

\proof
We proceed by double induction on $d$ and $d - r$.
If $d - r = 0$,
then $r = d$
so $\sumi$ is primary to $(x_1, \ldots, x_d)$.
In that case,
by assumption on $l$,
$(x_1^l, \ldots, x_d^l) \subseteq \sumi$.
Hence
$(x_1, \ldots, x_d)^{dln} \subseteq
(x_1^{nl}, \ldots, x_d^{ln})$.
As $x_j^l$ lies in some $I_p$,
it follows that
$x_j^{ln}$ lies in $I_p^n$,
hence
$(x_1, \ldots, x_d)^{dln} \subseteq \sumin$.
By Bayer-Stillman [BS, Lemma 1.7]
it follows now that
$\reg \sumipn \le dln$,
which proves the theorem in the case $d = r$.
Note that this also proves the case $d = 1$.

Now assume that $d > 1$, $r < d$.
By Lemma~\lmlm,
$$
\reg \sumipn
\le
\max \left\{
\reg \left(
\sumipn : x_{r+1}^{ln}
\right) +ln,
\reg \left(
\sumin + x_{r+1}^{ln} R
\right) \right\}.
$$
%HERE: I added Now
Now we use some facts about monomial ideals
(which definitely fail for arbitrary ideals).
First of all,
$\sumipn : x_{r+1}^{ln} R
= \sum_{j=1}^m (I_j^n : x_{r+1}^{ln})$.
Moreover,
by the choice of $l$,
$I_j^n : x_{r+1}^{ln} =
(I_j : x_{r+1}^{l})^n$
can be identified with the ideal $I_j^n$ after we set $x_{r+1}$ to $1$.

Thus to find the regularity of 
$\sum_j (I_j : x_{r+1}^{l})^n$,
we may work in the polynomial ring
$k[x_1, \ldots, x_r, x_{r+2}, \ldots, x_d]$
with one fewer variable.
Let $I = \sum_{j=1}^m I_{j}$ and
set
$$
\eqalign{
L_1 &=
\max_{S \subseteq \{x_{r+2}, \ldots, x_d\}}
\max_{x_q \in S}
\{
(d-1-|S|-r)l + (m+|S|-3+d) \reg\left(
I :x_{r+1}^lx_q^l) + I_{S,x_q}
\right)
\}, \cr
L_2 &=
\max_{S \subseteq \{x_{r+2}, \ldots, x_d\}}
\max_{x_q \in S}
\{
(d-|S|-r-1)l + (m+|S|-1+d) \reg\left(
(I :x_q^l) + x_{r+1}^lR + I_{S,x_q}
\right)
\}. \cr
}
$$
Thus by induction on $d$,
$$
\reg \left(\sumipn : x_{r+1}^{ln}\right)
=
\reg \sum_{j=1}^m (I_j : x_{r+1}^l)^n
\le n \max \{(d-1)l, L_1 \}.
$$
and by induction on $d-r$,
$\reg (\sumin + x_{r+1}^{ln}R) \le n \max \{dl, L_2\}$.
Hence
$$
\reg \sumipn
\le n \max \{ dl, L_1 + l, L_2\}.
$$

%HERE I rephrased this paragraph
Note that $L_2$ equals
$$
\max_{x_{r+1} \in S \subseteq \{x_{r+1}, x_{r+2}, \ldots, x_d\}}
\max_{x_q \in S\setminus\{x_{r+1}\}}
\{
(d-|S|-r)l + (m+|S|-2+d) \reg\left(
(I :x_q^l) + I_{S,x_q}
\right)
\}.
$$
and this is bounded above by $L$.

Now we want to show that $L_1 + l \le \max\{dl, L\}$.
As $x_{r+1}$ is not an element of $S$,
$(I :x_{r+1}^lx_q^l) + I_{S,x_q}
=
((I :x_q^l) + I_{S,x_q}): x_{r+1}^l$,
so that by Lemma~\lmlm,
the regularity of $I :x_{r+1}^lx_q^l + I_{S,x_q}$
is bounded above by
$$
\max\{ 0,
\reg ((I :x_q^l) + I_{S,x_q}) -l,
\reg (((I :x_q^l) + I_{S,x_q})+ x_{r+1}^l)+1 -l \}.
$$
Thus $L_1 + l$ is bounded above by the maximum of
$(d-|S|-r)l$,
$$
(d-|S|-r)l + (m+|S|-3+d) (\reg ((I :x_q^l) + I_{S,x_q}) -l)
$$
and
$$
(d-|S|-r)l + (m+|S|-3+d) (\reg (((I :x_q^l) + I_{S,x_q})+ x_{r+1}^l)+1 -l),
$$
as $S$ varies over subsets of
$\{x_{r+2}, \ldots, x_d\}$
and $x_q$ over elements of $S$.
Each of these expressions is clearly bounded above by
$\max\{dl, L\}$,
which finishes the proof of the theorem.
\qed

\vskip 3em
{\bf Acknowledgements:}
The authors thank the NSF for partial support.

%\vfill\eject
\bigskip
\bigskip
\centerline{\bf Bibliography}
\bigskip
\baselineskip=9pt
{%\sevenpoint

\item{[BS]}
D.\ Bayer and M.\ Stillman,
A criterion for detecting m-regularity,
{\it Invent.\ Math.}, {\bf 87} (1987), {1-11}.

\item{[EH]}
J.\ Eagon and M.\ Hochster,
R-sequences and indeterminates,
{\it Quarterly J.\ Math.}, {\bf 25} (1974), {61-71}.

\item{[E]}
D.\ Eisenbud,
Commutative Algebra with a View toward Algebraic Geometry,
Springer-Verlag, 1994.

\item{[EK]}
S.\ Eliahou and M.\ Kervaire,
Minimal resolutions of some monomial ideals,
{\it J.\ Algebra}, {\bf 129} (1990), 1-25.

\item{[HS]}
W.\ Heinzer and I.\ Swanson,
Ideals contracted from 1-dimensional overrings with an application
to the primary decomposition of ideals,
preprint.

\item{[HRS]}
W.\ Heinzer, L.\ J.\ Ratliff, Jr.\ and K.\ Shah,
Parametric decomposition of monomial ideals (I),
preprint.

\item{[HH]}
C.\ Huneke, and M.\ Hochster,
Tight closure, invariant theory, and the Brian\c con-Skoda theorem,
{\it Jour. of Am. Math. Soc.} {\bf 3}
(1990), {31-116}.

\item{[K]}
M.\ Katzman,
Finiteness of $\cup_e \hbox{Ass}\,F^e(M)$
and its connections to tight closure,
{\it Illinois J.\ Math.}, {\bf 40} (1996), no.2, 330-337.

\item{[K2]}
M.\ Katzman,
An application of Gr\"obner bases to a tight closure question,
preprint.

\item{[R]}
L.\ J.\ Ratliff, Jr.,
On prime divisors of $I^n$, $n$ large,
{\it Michigan Math.\ J.}, {\bf 23} (1976), {337-352}.

\item{[STV]}
B.\ Sturmfels, N.\ V.\ Trung and W.\ Vogel,
Bounds on degrees of projective schemes,
{\it Math.\ Annalen}, {\bf 302} (1995), {417-432}.

\item{[S]}
I.\ Swanson,
Powers of Ideals:
Primary decompositions, Artin-Rees lemma and regularity,
{\it Math.\ Annalen}, to appear.

\item{[T]}
W. \ Traves,
Differential operators on monomial rings,
preprint.

%\bigskip
%\bigskip
%\baselineskip=9pt
%\hbox{
%\vbox{\halign{#\hfil\cr
%Department of Mathematics \cr
%Massachusetts Institute of Technology 2-167 \cr
%Cambridge, MA 02139 \cr
%{\tt kesmith@math.mit.edu} \cr
%}}
%\hskip 1.1in
%\vbox{\halign{#\hfil\cr
%Department of Mathematical Sciences \cr
%New Mexico State University \cr
%Las Cruces, NM 88003-8001 \cr
%{\tt iswanson@nmsu.edu} \cr
%}}
%}

\end